\documentclass[12pt]{amsart}
\usepackage[]{amsmath, amsthm, amsfonts, amssymb, epsfig}
\newcommand\C{{\mathbb C}}

\newcommand\dee{\partial}
\newcommand\Om{\Omega}

\renewcommand\phi{\varphi}

\numberwithin{equation}{section}

\begin{document}

\title[Poisson and Dirichlet]{Something about Poisson and Dirichlet}
\author[S.~R.~Bell and L.~R.~de~la~Torre]
{Steven R.~Bell and Luis Reyna de la Torre}

\address[]{Mathematics Department, Purdue University, West Lafayette,
IN 47907}
\email{bell@math.purdue.edu}

\address[]{Department of Physics, Arizona State University, Tempe,
AZ 85281}
\email{ lereyna@asu.edu}



\begin{abstract}
We solve the Dirichlet problem in the unit
disc and derive the Poisson formula using
very elementary methods and explore consequent
simplifications in other foundational areas
of complex analysis.
\end{abstract}

\maketitle

\theoremstyle{plain}

\newtheorem {thm}{Theorem}[section]
\newtheorem {lem}[thm]{Lemma}

\hyphenation{bi-hol-o-mor-phic}
\hyphenation{hol-o-mor-phic}

\section{Mathematical DNA}
\label{sec1}

For reasons that have always been mysterious to the first
author, he has often found himself thinking about the
Dirichlet problem in the plane and has fought urges to
look for explicit formulas for the Poisson kernels
associated to various kinds of multiply connected domains,
especially quadrature domains. He has obsessed about
solutions to the Dirichlet problem with rational boundary
data (see \cite{BEKS}) and he cannot stop thinking about the
Khavinson-Shapiro conjecture about the same problem with
polynomial data. His publication list is interspersed with
papers where he has given in and found formulas for the
Poisson kernel in terms of the Szeg\H o kernel and in
terms of Ahlfors maps. (See
\cite{{Be,Be2}} for an expository treatment of some of these
results.) After thinking about the Dirichlet problem his
entire adult life, he can solve it a different way every
day of the week. He recently looked up his mathematical
lineage at the Math Genealogy Project and found a possible
explanation for his obsession. He is a direct mathematical
descendent of both Poisson and Dirichlet. These problems
are in his blood!

The authors worked together on a summer research project at
Purdue University in 2018 to find a particularly elegant
and simple way to approach these problems. We want to demonstrate
here how Poisson and Dirichlet might solve their famous
problems today if they had lived another 200 years and
developed a major lazy streak. We assume
that our reader has seen a traditional approach to this subject
in a course on complex analysis and so will appreciate the
novelty and smooth sailing of the line of reasoning here, but
just in case the reader hasn't, we have tried to present the
material in a way that can be understood assuming only
a background in basic analysis.

We would like to thank Harold Boas for reading an early
draft of this work and making many valuable suggestions
for improvement. Harold knows a lot about the family
business because he, like the first author, was a student
of Norberto Kerzman at MIT in the late 1970's.

\section{Harmonic functions}
\label{sec2}

We define a harmonic function on a domain
$\Om$ in the complex plane to be a continuous complex
valued function $u(z)$ on the domain that satisfies the
{\it averaging property\/} on $\Om$, meaning that
$$u(a)=\frac{1}{2\pi}\int_0^{2\pi} u(a+re^{i\theta})\ d\theta$$
whenever $D_r(a)$, the disc of radius $r$ about $a$, is
compactly contained in $\Om$. It is well known that
this definition of harmonic function is equivalent to
all the other standard definitions (see, for example,
Rudin \cite[Chap.~11]{Ru}), and we will demonstrate
this in very short order. Since we will
explore other definitions of harmonic, we will
emphasize that we are currently thinking of harmonic
functions as being defined in terms of an averaging
property by calling them harmonic-ave functions.

We first note that analytic polynomials are harmonic-ave.
Indeed, if $P(z)=a_n z^n+a_{n-1}z^{n-1}+\dots +a_1 z+a_0$,
the averaging property is clear on discs centered at the
origin because the constant function $a_0$ obviously
satisfies the averaging property at the origin and so
does $z^n$ for $n\ge1$ because
$$\int_0^{2\pi} (re^{i\theta})^n\ d\theta
=r^n\int_0^{2\pi} \cos(n\theta) \ d\theta
+ir^n\int_0^{2\pi} \sin(n\theta) \ d\theta =0=0^n.$$
To see that the averaging property holds at a point
$a\in\C$, write
$$P(z)=P((z-a)+a)$$
and expand to get a polynomial in $(z-a)$. Now the
argument we used at the origin can be applied to
discs centered at $a$.

It follows that the averaging property also holds
for conjugates of polynomials in $z$, so polynomials
in $\bar z$ are harmonic-ave. Note also that complex
valued functions are harmonic-ave if and only if
their real and imaginary parts are both harmonic-ave.

The {\it Dirichlet problem\/} on the unit disc is:
given a continuous real valued function $\phi$ on
the unit circle, find a real valued function $u$
that is continuous on the closure of the unit
disc with boundary values given by $\phi$ such that
$u$ is harmonic on $D_1(0)$. We will find a
solution to this problem that is harmonic in the
averaging sense in the next section. To do so, we
will need to know the elementary fact that a
continuous real valued function on the unit circle
can be uniformly approximated on the unit circle
by a real polynomial $p(x,y)$. We will now
demonstrate this little fact, assuming the
Weierstrass theorem about the density of real
polynomials of one variable among continuous
functions on closed subintervals of the real line.

Suppose we are given a continuous real valued
function $\phi$ on the unit circle. We wish
to find a real polynomial $p(x,y)$ that is
uniformly close to $\phi$ on the unit circle.
Note that we may assume that $\phi(\pm 1)=0$ 
because we may subtract a polynomial function
of the form $ax+b$ to make $\phi$ zero at $\pm 1$.
Define two continuous functions on $[-1,1]$ via
$$h_{\text{top}}(x)=\phi(x + i\sqrt{1-x^2})$$ 
and
$$h_{\text{bot}}(x)=\phi(x - i\sqrt{1-x^2}).$$
We can uniformly approximate
$h_{\text{top}}$
and
$h_{\text{bot}}$
on $[-1,1]$ by real polynomials
$p_{\text{top}}(x)$
and
$p_{\text{bot}}(x)$.
Next, let $\chi_\epsilon(y)$ be a continuous
function that is equal to one for $y>\epsilon$,
equal to zero for $y<-\epsilon$ and follows
the line connecting $(-\epsilon,0)$ to $(\epsilon,1)$
for $-\epsilon\le y\le \epsilon$. Let $p_\epsilon(y)$
be a polynomial in $y$ that is uniformly close to
$\chi_\epsilon$ on $[-1,1]$. Now, because $\phi$
vanishes at $\pm 1$, the polynomial
$$p_{\text{top}}(x) p_\epsilon(y)+
p_{\text{bot}}(x) p_\epsilon(-y)$$
approximates $\phi$ on the unit circle, and the
approximation can be improved uniformly by
shrinking $\epsilon$ and improving the approximations
of the other functions involved.

We are now in position to solve the Dirichlet
problem on the disc, but before we begin in
earnest, this is a good place to emphasize
that Green's theorem for the unit disc
depends on nothing more than the fundamental
theorem of calculus from freshman calculus.
Indeed, let $P(x,y)$ be a $C^1$-smooth function.
Let $C_{\text{top}}$ denote the top half of the
unit circle parameterized in the clockwise sense
by $z(x)=x+iy_{\text{top}}(x)$, $-1\le x\le 1$,
where $y_{\text{top}}(x)=\sqrt{1-x^2}$, and let
$C_{\text{bot}}$ denote the bottom half of the
unit circle parameterized in the counterclockwise
sense by $z(x)=x+iy_{\text{bot}}(x)$, $-1\le x\le 1$,
where $y_{\text{bot}}(x)=-\sqrt{1-x^2}$.
Note that the unit circle $C_1(0)$ parametrized
in the counterclockwise sense is given by 
$C_{\text{bot}}$ followed by $-C_{\text{top}}$.
Drum roll\dots
\begin{align*}
\int_{C_1(0)} P\,dx &=
\left(\int_{C_{\text{bot}}} -\int_{C_{\text{top}}}\right) P\,dx=
\int_{-1}^1 \left[P(x,y_{\text{bot}}(x))
- P(x,y_{\text{top}}(x))\right] \ dx \\
&=-
\int_{-1}^1\left(\int_{y_{\text{bot}}(x)}^
{y_{\text{top}}(x)}\frac{\dee P}{\dee y}(x,y)\ dy\right) dx
=\iint_{D_1(0)} -\frac{\dee P}{\dee y} \ dx\wedge dy.
\end{align*}
The other half of Green's formula follows by repeating
the argument using the words left and right in place of
top and bottom. This argument on the disc can be easily
generalized to demonstrate Green's theorem on any
region that can be cut up into regions that have a top
boundary curve and a bottom curve and a left curve and
a right curve. An annulus centered at the origin cut
into four regions by the two coordinate axes is such
a domain. We will need Green's theorem later in the paper
when we study analytic functions from a philosophical point
of view inspired by our observations about harmonic
functions and the Dirichlet problem.

\section{Solution of the Dirichlet problem on the unit disc}
\label{sec3}
The unit disc has the special feature that, given
polynomial data $\phi$, it is straightforward to
write down a polynomial solution to the Dirichlet
problem with boundary values given by $\phi$.
Indeed, a polynomial $p(x,y)$ in the real variables
$x$ and $y$ can be converted to a polynomial in $z$
and $\bar z$ by replacing $x$ by $(z + \bar z)/2$ and
$y$ by $(z - \bar z)/(2i)$ and expanding. It is now
an easy matter to extend the individual
terms in the sum to harmonic-ave functions on the
disc by noting that
$$z^n{\bar z}^{\,m}$$
is equal to one on the unit circle if $n=m$, equal to
$z^{n-m}$ on the circle if $n>m$, and equal to
${\bar z}^{\,m-n}$ on the circle if $m>n$, each of which
is harmonic-ave inside the unit circle.

Now, given a continuous real valued function $\phi$
on the unit circle, there is a sequence of real
valued polynomials $p_n(x,y)$ that converges uniformly
to $\phi$ on the unit circle. Let $u_n$ be the polynomial
harmonic-ave extension of $p_n$ to the disc described in the
paragraph above. Note that $u_n$ can be expressed as a
constant plus a polynomial in $z$ that vanishes at the
origin plus a polynomial in $\bar z$ that also vanishes
at the origin. We now claim that the functions $u_n$
converge uniformly on the closed disc to a solution of
the Dirichlet problem. To see this, we must first show
that real valued harmonic-ave functions $u$ on the disc
that extend continuously to the closure satisfy the
{\it maximum principle\/} in the form
$$\max \{u(z): |z|\le 1\} =
 \max \{u(e^{i\theta}):0\le \theta\le 2\pi\}.$$
Indeed, if the maximum value of such a function $u$
occurs at a point $z_0$ inside the unit circle, we
can express the value of $u$ at $z_0$ as an average
of $u$ over a small circle centered at $z_0$. We can
let the radius of that circle increase until the
circle touches the unit circle at a single point.
The averaging property holds on the limiting circle
because of uniform continuity. Let $M$ denote the
maximum value $u(z_0)$. Now, in order for the average
of a continuous function that is less than
or equal to $M$ over that circle to be equal to $M$,
it must be that $u$ is equal to $M$ on the whole
circle. Hence the value of $u$ at the point where the
inner circle touches the unit disc must also be $M$.
This proves the maximum principle inequality for real
harmonic-ave functions. The minimum principle follows by
applying the maximum principle to $-u$.

Because the $p_n$ converge to $\phi$ uniformly on the
unit circle, the sequence $\{p_n\}$ is uniformly Cauchy
on the unit circle, i.e., given $\epsilon>0$, there is
an $N$ such that $|p_n-p_m|<\epsilon$ on the unit
circle when $n$ and $m$ are greater than $N$. The
maximum and minimum principle inequalities applied
to the imaginary parts of $u_n$ (which are zero on
the unit circle) allow us to conclude that the
functions $u_n$ are {\it real\/} valued.
Furthermore, the maximum and minimum principles
applied to $u_n-u_m$ show that the uniformly Cauchy
estimates for the sequence $\{p_n\}$ on the unit
circle extend to hold for the sequence $\{u_n\}$ on
the whole closed unit disc, showing that $\{u_n\}$
is uniformly Cauchy on the closed unit disc. Hence,
since the $u_n$ are continuous, they converge
uniformly on the closed unit disc to a continuous
function $u$ that is equal to $\phi$ on the unit circle.
Finally, it is clear that $u$ is harmonic-ave on the
inside of the unit circle because the averaging
property is preserved under uniform limits. We have
solved the Dirichlet problem in a purely existential
manner without ever differentiating a function! We
now turn to Poisson's problem of finding a
{\it formula\/} for our solution.

\section{Poisson's formula}
\label{sec4}

Notice that the set of harmonic-ave functions
$$\{1, z^n, {\bar z}^{\,n} : n=1,2,3, \dots\}$$
is orthonormal under the inner product
$$\langle u,v\rangle =
\frac{1}{2\pi}\int_0^{2\pi} u(e^{i\theta})\ 
\overline{v(e^{i\theta})}\ d\theta.$$
Define $K_N(z,w)$ via
$$K_N(z,w):= 1 + \sum_{n=1}^N z^n {\bar w}^{\,n} +
\sum_{n=1}^N {\bar z}^{\,n} w^n,$$
and observe that if $u(z)=z^n$, then
\begin{equation}
\label{poissonN}
u(z)=
\frac{1}{2\pi}\int_0^{2\pi}
K_N(z,e^{i\theta})\, u(e^{i\theta}) \ d\theta
\end{equation}
for $z\in D_1(0)$ if $N\ge n$ because of the
orthonormality of the terms in the sum. The
same is true if $u(z)\equiv 1$ or $u(z)={\bar z}^{\,n}$.
We conclude that, if $u$ is the solution of
the Dirichlet problem for polynomial data
$p$ of degree $n$ as constructed in the previous
section, then formula (\ref{poissonN}) holds
for $u(z)$ for $z\in D_1(0)$ when $N>n$.

Using the famous geometric series formula,
$$
1+\zeta+\dots+\zeta^N
=
\frac{1}{1-\zeta}- \frac{\zeta^{N+1}}{1-\zeta},
$$
we see that
$$K_N(z,w)= 1 + (z\,\bar w)
\sum_{n=0}^{N-1} z^n {\bar w}^{\,n} +
 (\bar z\, w)
\sum_{n=0}^{N-1} {\bar z}^{\,n} w^n$$
converges uniformly in
$w=e^{i\theta}$ when $z\in D_1(0)$ to
$$K(z,w):=1+\frac{z\,\bar w}{1-z\,\bar w}+
\frac{{\bar z}\,w}{1-{\bar z}\,w},$$
and the error $\mathcal E_N=|K_N-K|$ is controlled via
$$\mathcal E_N(z,w)\le \frac{2|z|^{N+1}}{1-|z|}$$
when $z\in D_1(0)$ and $|w|=1$.
Hence, it follows by taking uniform limits that
formula (\ref{poissonN}) holds for $z\in D_1(0)$
with the $N$ removed for the polynomials $u_n$ that
we constructed from polynomial boundary data $p_n$.
We can now let the polynomials $p_n$ tend uniformly
to $\phi$ and use the fact proved above that the
corresponding solutions $u_n$ to the Dirichlet
problem with boundary data $p_n$ converge uniformly
to a solution $u$ of the Dirichlet problem to obtain
Poisson's famous formula for the solution to the
Dirichlet problem,
\begin{equation}
\label{poisson}
u(z)=
\frac{1}{2\pi}\int_0^{2\pi}
K(z,e^{i\theta}) \phi(e^{i\theta}) \ d\theta.
\end{equation}
This formula reveals that the solution $u$ can
be written
$$u(z)=a_0+h(z)+\overline{h(z)}$$
where $a_0$ is the (real valued) average of $\phi$ on the
unit circle and $h$ is an analytic function on $D_1(0)$
that vanishes at the origin given via
$$h(z)=\frac{1}{2\pi}\int_0^{2\pi}
\frac{z e^{-i\theta}}{1-z e^{-i\theta}}\ \phi(e^{i\theta})\ d\theta
=
\frac{z}{2\pi i}\int_{C_1}
\frac{\phi(w)}{w(w-z)} \ dw,$$
where $C_1$ denotes the unit circle parameterized in
the standard sense using $w=e^{i\theta}$ and
$dw=ie^{i\theta}\,d\theta$. It is now a rather
easy exercise to take limits of complex difference
quotients to see that complex derivatives in
$z$ can be taken under the integral sign in the
definition of $h$. Hence, $h(z)$ is infinitely
complex differentiable. It follows from the
Cauchy-Riemann equations that $u$ is a
$C^\infty$-smooth real valued function on
$D_1(0)$ in $x$ and $y$ that satisfies the Laplace
equation there and solves the Dirichlet problem.
Furthermore, $u$ is the real part of an infinitely
complex differentiable function $H=a_0+h/2$
on $D_1(0)$.

The maximum principle yields that $u$ is the
{\it unique\/} solution to the problem in the
realm of harmonic functions understood in the
sense of averaging. We will see in the next
section that it is also the unique solution
among harmonic functions defined in the traditional
sense of satisfying the Laplace equation.

We remark here that rather simple algebra reveals the
well-known formulas for the Poisson kernel,
$$K(z,w) =\text{Re }\frac{1+z\,\bar w}{1-z\,\bar w}
=\text{Re }\frac{w+z}{w-z}$$
and
$$K(z,e^{i\theta}) =\frac{1-|z|^2}{|z-e^{i\theta}|^2}.$$

It is a routine matter to extend the above
line of reasoning to any disc (either by repeating
the argument or making a complex linear change
of variables $Az+B$). Since a complex valued
function is harmonic-ave if and only if its real
and imaginary parts are harmonic-ave, it follows
from our work that a complex valued harmonic-ave
function is given locally by $g+\overline{G}$ where
$g$ and $G$ are infinitely complex differentiable
functions.

At this point, it would be tempting to experiment
with thinking of analytic functions as being
harmonic-ave functions that do not involve the
antianalytic $\overline{G}$ parts. The formula for
$h$ above reveals that the analytic $g$ part is
locally a uniform limit of analytic polynomials.
Since $\{1, z^n: n=1,2,3,\dots\}$ are
orthonormal on the unit circle, we could let
$$k_N(z,w)=1+\sum_{n=1}^N z^n{\bar w}^{\,n}$$
and use the same line of reasoning that we used earlier
in this section to conclude that
\begin{equation}
\label{CauchyN}
f(z)=\frac{1}{2\pi}
\int_0^{2\pi}k_N(z,e^{i\theta})\,f(e^{i\theta})\ d\theta
\end{equation}
if $f(z)$ is equal to $1$ or $z^n$ with $n\le N$.
The geometric series estimate we used above
shows that $k_N(z,w)$ converges uniformly in $w$
on the unit circle for fixed $z$ in $D_1(0)$ to
$$k(z,w) := \frac{1}{1-z\,\bar w}.$$
Hence, we can let $N\to\infty$ in (\ref{CauchyN})
to see that
\begin{equation}
\label{Cauchy}
f(z)=\frac{1}{2\pi}
\int_0^{2\pi}k(z,e^{i\theta})\,f(e^{i\theta})\ d\theta
\end{equation}
when $f$ is an analytic polynomial. Finally, if
$f$ is a uniform limit of analytic polynomials
on an open set containing the closed unit disc,
we may conclude that $f$ satisfies identity
(\ref{Cauchy}), too. The identity can easily be seen
to be the classical Cauchy integral formula
on the unit disc, and from this point, the theory of
analytic functions would gush forth. In particular,
analytic functions would be seen to be infinitely complex
differentiable and given locally by convergent power
series. We will explore this idea and various other
alternate ways of thinking about analytic functions
after we verify that our definition of harmonic
functions via an averaging property gives rise to
the same set of functions as any of the more
standard definitions.

\section{Traditional definitions of harmonic functions}
\label{sec5}

Some complex analysis books define harmonic functions
to be twice continuously differentiable functions
that satisfy the Laplace equation. With this definition,
one can use Green's identities on an annulus (which we
pointed out in section~1 to be quite elementary) to show
that such harmonic functions satisfy the averaging property.
Hence, this class of functions could be seen to be the
same as harmonic-ave functions. We won't pursue this
idea here because we can easily prove something stronger
with less effort.

It is most gratifying to define harmonic functions to
be merely continuous functions $u$ whose first partial
derivatives exist and whose second partial derivatives
$\dee^2 u/\dee x^2$ and $\dee^2 u/\dee y^2$ exist and
satisfy the Laplace equation. Call such functions
harmonic-pde. We will now adapt a classic
argument to show that the class of harmonic-pde
functions agrees with our class of continuous functions
that satisfy the averaging property on an open set. Indeed,
if $u$ is a harmonic-pde function defined on an open set
containing the closed unit disc, then $u$ minus the
harmonic-ave function $U$ that we constructed solving
the Dirichlet problem on the unit disc with the same
boundary values as $u$ on the unit circle, if not the
zero function, would have either a positive maximum or
a negative minimum. Suppose it has a positive maximum
$M>0$ at a point $z_0$ in $D_1(0)$. Choose $\epsilon$
with $0<\epsilon< M$. Now
$$v(z) := u(z)-U(z) + \epsilon |z|^2$$
is equal to $\epsilon$ on the unit circle and
attains a positive value at $z_0$ that is larger
than $\epsilon$. Hence, $v$
attains a positive maximum at some point
$w_0$ in $D_1(0)$. The Laplacian of $v$ at $w_0$
is equal to $4\epsilon$, which is strictly positive.
However, the one variable second derivative test
from freshman calculus applied to $v$ in the $x$-direction
yields that $\dee^2 v/\dee x^2$ must be less
than or equal to zero. (If it were positive,
$v$ could not have a local maximum at $w_0$.)
Similarly, the second derivative test in the
$y$-direction yields that $\dee^2 v/\dee y^2$
must be less than or equal to zero. We conclude
that the Laplacian of $v$ at $w_0$ must be
less than or equal to zero, which is a contradiction
because the Laplacian is strictly positive on
$D_1(0)$. This shows that $u-U$ cannot
have a positive value. If we replace
$u-U$ by $U-u$, the same reasoning shows
that $U-u$ cannot have a positive value.
Hence $U-u\equiv 0$ on the unit disc and
we have shown that $u$, like $U$, is harmonic in
the sense of being averaging. Finally,
because the operations of translating
and scaling preserve harmonic functions in
both senses of the word, we can translate and
scale any disc to the unit disc to be able to
deduce the equivalence of the two definitions
of harmonic on any open set. We now stop
hyphenating harmonic and turn to hyphenating
analytic.

\section{Analytic functions}
\label{sec6}
When the 200+ year old Poisson and Dirichlet
took a break from harmonic functions, they might
have turned their attention to similar considerations
applied to analytic functions, which satisfy
both the averaging property on discs and Cauchy's
theorem for circles.

In order to understand the rest of this section,
the reader will need to know (or look up) the
definition of the complex contour integral
$\int_\gamma f\,dz$ along a curve $\gamma$, the
basic estimate
$$\left|\int_\gamma f\,dz\right|\le
\sup \{|f(z)|: z\in\gamma\}\cdot \text{ Length}(\gamma),$$
and the fundamental theorem of calculus
$\int_\gamma f'\,dz=f(b)-f(a)$
for complex contour integrals, assuming that $f$ is
continuously complex differentiable and that $\gamma$
starts at $a$ and ends at $b$. We will also use
the differential operators
$$\frac{\dee}{\dee z} = \frac{1}{2}\left(
\frac{\dee}{\dee x} - i
\frac{\dee}{\dee y}
\right)
\qquad
\text{ and }
\qquad
\frac{\dee}{\dee\bar z} = \frac{1}{2}\left(
\frac{\dee}{\dee x}+ i
\frac{\dee}{\dee y}
\right).$$
These two very important operators can be ``discovered''
by writing $dz=dx+i\,dy$ and $d\bar z=dx-i\,dy$ and
manipulating
$$df = \frac{\dee f}{\dee x}\,dx +
\frac{\dee f}{\dee y}\,dy$$
to appear in the form
$$df = \frac{\dee f}{\dee z}\,dz +
 \frac{\dee f}{\dee\bar z}\,d\bar z.$$
The condition $\frac{\dee f}{\dee\bar z}=0$ is equivalent
to writing the Cauchy-Riemann equations for the real
and imaginary parts of $f$. If $h(z)$ is a complex
differentiable function, then
$\frac{\dee h}{\dee\bar z}=0$ and
$\frac{\dee h}{\dee z}=h'(z)$. Furthermore,
$\frac{\dee \overline{h}}{\dee z}=0$ and
$\frac{\dee \overline{h}}{\dee\bar z}=\overline{h'(z)}$.

We will denote the boundary circle of the
disc $D_r(a)$ parameterized in the counterclockwise
sense by $C_r(a)$. Writing out the real and imaginary
parts of the contour integral of a complex function
$f$ around $C_r(a)$ and applying Green's
theorem for a disc to the real and imaginary parts
yields, what we like to call, the {\it complex Green's
theorem\/} for a disc,
$$\int_{C_r(z_0)} f\ dz=
2i\iint_{D_r(z_0)} \frac{\dee f}{\dee\bar z}\ dA,$$
where $dA$ denotes the element of area $dx\wedge dy$.
(Using the real Green's theorem to prove the complex
Green's theorem is quick and easy, but a more
civilized way to deduce the theorem would be to note
that
$$d(f\,dz)=
\frac{df}{dz}\ dz\wedge dz+
\frac{df}{d\bar z}\ d\bar z\wedge dz=
0+\frac{df}{d\bar z}\,(2i dx\wedge dy)$$
and to apply Stokes' theorem.)

Define a continuous complex valued function $f$ to be
analytic-circ on a domain $\Omega$ if it
is harmonic on $\Om$ (and so satisfies the averaging
property) and, for each closed subdisc of $\Omega$,
the complex contour integral of $f$ over the boundary
circle is zero.

Complex polynomials are easily seen to be
analytic-circ because the monomials are complex
derivatives of monomials one degree higher and the
fundamental theorem of calculus for complex contour
integrals reveals that the integrals around closed
curves are zero. It now follows that functions that
are the uniform limit of complex polynomials
(in $z$) on each closed subdisc of $\Omega$
are analytic. We will now show that functions
that are analytic-circ must be the uniform limit
of complex polynomials on each closed subdisc of
$\Omega$.

Suppose $f$ is analytic-circ on $\Om$. Let $D_r(z_0)$
be a disc that is compactly contained in $\Om$ and
let $C_r(z_0)$ denote the boundary circle
parameterized in the counterclockwise sense. The
complex Green's theorem for a disc yields that
$$0=\int_{C_r(z_0)} f\ dz=
2i\iint_{D_r(z_0)} \frac{\dee f}{\dee\bar z}\ dA,$$
where $dA$ denotes the element of area $dx\wedge dy$.
Since this integral is zero for every disc compactly
contained in $\Om$, it follows that 
$\dee f/\dee\bar z\equiv0$ on $\Om$, i.e., that the
real and imaginary parts of $f$ satisfy the Cauchy-Riemann
equations. Since harmonic functions are infinitely
differentiable, the textbook proof that $f$ is
complex differentiable can now be applied.

Since being analytic-circ is a local property
and is invariant under changes of variables of
the form $Az+B$, we may restrict our attention
to a function $f$ that is analytic-circ
on a neighborhood of the closed unit disc.
We can gain more insight into the implications
of the definition by noting that such a harmonic
function is given by $g+\overline{G}$ where
$g$ and $G$ are infinitely complex
differentiable. It follows from
our complex Green's calculation above that
$\dee f/\dee\bar z=\overline{G'}\equiv0$
on $\Om$, and that $f$ is therefore equal to
the analytic part $g$ plus a constant. From
this point, it follows that $f$ is given by a
Cauchy integral formula and we merge into the
fast lane of the classical theory of analytic
functions.

Of course, the traditional way to define analytic
functions is as complex differentiable functions on
open sets. Call such functions analytic-diff. We
will now show that the traditional definition leads
to the same class of functions that we have defined
as being analytic-circ.

Suppose that $f$ is an analytic-diff function
defined on an open set containing the closed unit
square $\mathcal S:=[0,1]\times[0,1]$.
We will show that the complex integral of $f$ around
the counterclockwise perimeter curve $\sigma$ of the
square must be zero. The well known argument will be
a beautiful bisection method tracing back to Goursat.

It follows from our assumption that $f$ can be
locally well approximated by a complex linear function
in the following sense. Suppose $a$ is a point in
$\mathcal S$, and let $\epsilon >0$ be given.
Since $f'(a)$ exists, there is a $\delta>0$ such that
$$\frac{f(z)-f(a)}{z-a}=f'(a) + E_a(z)$$
where $|E_a(z)|<\epsilon$ when $|z-a|<\delta$, $z\ne a$.
Define $E_a(a)$ to be zero to make $E_a$ continuous
at $a$ and so as to be able to assert that
$$f(z)=f(a)+f'(a)(z-a) + E_a(z)(z-a)$$
on the whole open set where $f$ is defined and $E_a$ is
a continuous function in $z$ on that set. Furthermore, 
$|E_a(z)|<\epsilon$ on $D_\delta(a)$.
The complex integral of the polynomial $f(a)+f'(a)(z-a)$
around any square is zero because first degree
polynomials are derivatives of second degree
polynomials. If $\sigma_h$ is the counterclockwise
boundary curve of a small square ${\mathcal S}_h$ of
side $h$ contained in $D_\delta(a)$ that
contains the point $a$, it follows that
$$\left|\int_{\sigma_h}f(z)\ dz\right|=
\left|\int_{\sigma_h}E_a(z)(z-a)\ dz\right|
\le \epsilon (\sqrt{2}h)(4h).$$
Hence,
\begin{equation}
\label{toosmall}
\left|\int_{\sigma_h} f\ dz\right|
\le
(4\sqrt{2}\,) \epsilon\,
\text{Area\,}({\mathcal S}_h).
\end{equation}
We will now follow a version of Goursat's famous
argument to explain how this could be made {\it too
small\/} if $\int_\sigma f\ dz$ were not zero.

Indeed, suppose that $I:=\int_\sigma f\ dz$ is not
equal to zero. Note that $I$ is equal to the sum
of the integrals around the four counterclockwise
squares obtained by cutting the big square into four
equal squares of side $1/2$ since the integrals
along the common edges cancel. For these four integrals
to add up to the non-zero value $I$, the modulus of at
least one of them must be greater than or equal to
$|I|/4$. Name such a square ${\mathcal S}_1$ and its
counterclockwise boundary curve $\sigma_1$. Note that
$$\left|\int_{\sigma_1} f\ dz\right|
\ge |I|\cdot\text{Area\,}({\mathcal S}_1).$$

We may now dice up ${\mathcal S}_1$ into four equal
subsquares and repeat the argument to obtain a square
${\mathcal S}_2$ with boundary curve $\sigma_2$
such that the modulus of the integral of $f$ around
$\sigma_2$ is greater than or equal to $|I|$ times the
area of ${\mathcal S}_2$. Continuing in this
manner, we obtain a nested sequence of closed
squares $\{{\mathcal S}_n\}_{n=1}^\infty$
with boundary curves $\sigma_n$, the diameters
of which tend to zero as $n\to\infty$, such that
\begin{equation}
\label{toobig}
\left|\int_{\sigma_n} f\ dz\right|
\ge |I|\cdot\text{Area\,}({\mathcal S}_n).
\end{equation}
There is a unique point $a$ that belongs to
all the squares. Now, given an $\epsilon$ less
than $|I|/(4\sqrt{2})$, the squares that
eventually fall in the disc $D_\delta(a)$ that
we specified above satisfy both area inequalities
(\ref{toosmall}) and (\ref{toobig}), which are
incompatible. This contradiction shows that $I$
must be zero!

Since any square can be mapped to the unit
square via a mapping of the form $Az+B$, it
follows from a simple change of variables that
the integral of an analytic-diff function
around any square must be zero. Furthermore,
any rectangle can be approximated by a rectangle
subdivided into a union of $n\times m$ squares.
It follows that the integral of an
analytic-diff function around any rectangle must
be zero.

From this point, there are several standard arguments
to prove the Cauchy integral formula on a disc for
such functions (see Ahlfors \cite[p.~109]{Ah} or
Stein \cite[p.~37]{St}). It then follows from the
Cauchy integral formula that such a function would be
locally the uniform limit of analytic polynomials,
and so the function would be analytic-circ. However,
we can simplify these standard arguments by using
some of the power of our work on harmonic functions
in the previous sections.

Suppose that $f$ is analytic-diff on an open disc.
Since $f$ is complex differentiable, it is continuous.
Define $F(z)$ at a point $z$ in the disc by the
integral of $f$ along a horizontal ``zig'' from the
center followed by a vertical ``zag'' connecting to the
point $z$. The fundamental theorem of calculus
reveals that
$$\frac{\dee}{\dee y} F(x+iy)=if(x+iy).$$
Since the integral of $f$ around rectangles is
zero, we could also define $F$ via an integral
along a vertical zag followed by a horizontal zig.
Using this definition, the fundamental theorem of
calculus shows that
$$\frac{\dee}{\dee x} F(x+iy)=f(x+iy).$$
Hence, $F$ is a continuously differentiable
function whose real and imaginary parts satisfy
the Cauchy-Riemann equations, and furthermore,
is a complex differentiable function such that
$F'=f$ on the disc. Repeat this construction to
get a twice continuously complex differentiable
function $G$ such that $G''=f$.

Now, since $G$ is twice continuously complex
differentiable, it is easy to use the
Cauchy-Riemann equations to show that the real
and imaginary parts of $G$ are harmonic functions.
Indeed, if $G(x+iy)=u(x,y) + i v(x,y)$, then
the Cauchy-Riemann equations yield that
$$G'=u_x+iv_x = v_y-iu_y$$
and
$$G''=u_{xx}+iv_{xx} = -u_{yy}-iv_{yy}$$
and we see that $u$ and $v$ satisfy the Laplace
equation by equating the real and imaginary parts
of $G''$. Our work in previous sections shows that
these harmonic functions are $C^\infty$-smooth and it
follows that the real and imaginary parts
of $f$ are $C^\infty$-smooth and satisfy the
Laplace equation and the Cauchy-Riemann equations.
We may now use Green's theorem on a disc to prove
the Cauchy theorem for $f$ on discs. Hence $f$
is analytic-circ and we have proved the equivalence
of the definitions, revealing that $f$ is locally
the uniform limit of complex polynomials and is
given by the Cauchy integral formula. The
shortcuts we have revealed in the theory of
analytic functions deliver us to page~114 of
Ahlfors.

Before we conclude this section, we present
one last alternate way to define analytic functions
that might be of interest to experienced analysts.
The result is known (see for example,
Springer \cite{Sp} or Globevnik \cite{Gl}), but we
are in a position to prove it rather efficiently
here. We now will define a continuous complex
valued function to be analytic-ave on the unit disc
if $f(z)(z-a)$ satisfies the averaging property on
circles $C_r(a)$ contained in the disc, i.e., if
$$0=\int_0^{2\pi} f(a+re^{i\theta})(re^{i\theta})\ d\theta$$
whenever the closure of $D_r(a)$ is contained in
the unit disc. Note that this condition is
equivalent to the condition that
$$0=\int_{C_r(a)} f\ dz$$
for each such circle. We will now prove that analytic-ave
functions are analytic in the usual sense. This result can
be viewed as a version of Morera's theorem saying that
a continuous complex valued function that satisfies the
Cauchy theorem on circles must be analytic.

Let $\chi(t)$ be an real valued non-negative function
in $C^\infty[0,1]$ that is equal to one for $t<\frac12$
and equal to zero for $t>\frac34$. Define
$$\phi(z)=c\chi(|z|^2)$$
where $c$ is chosen so that
$\int \phi \,dA=1$. Let $\phi_\epsilon$
denote the approximation to the identity given by
$$\phi_\epsilon(z)=\frac{1}{\epsilon^2}\phi(z/\epsilon).$$
The proof of our claim rests on a straightforward
calculation that shows that
$$\frac{\dee}{\dee\bar z}\phi_\epsilon(z-w)$$
is $(z-w)$ times a function $\psi_\epsilon(z-w)$
that is radially symmetric about $z$ in $w$.  In
fact,
$$\psi_\epsilon(z)=
\frac{c}{\epsilon^4}\chi'(|z|^2/\epsilon^2).$$
The calculation hinges on the chain rule plus the
fact that $$\frac{\dee }{\dee\bar z}|z^2|=
\frac{\dee }{\dee\bar z}(z\,\bar z)=z.$$

Given a continuous function $f$ on the unit disc
such that $f(z)(z-a)$ satisfies the averaging
property on circles $C_r(a)$ compactly contained
in the disc, let $f_\epsilon=\phi_\epsilon* f$
for small $\epsilon>0$. Note that $f_\epsilon$
is $C^\infty$ smooth on $D_{1-\epsilon}(0)$ and
that $f_\epsilon$ converges uniformly on compact
subsets of the unit disc to $f$ as $\epsilon\to0$.
One can differentiate under the integral in the
convolution formula to see that
$$\frac{\dee f_\epsilon}{\dee\bar z}=
\frac{\dee \phi_\epsilon}{\dee\bar z}*f
=\iint_{w\in D_1(0)}
\frac{\dee \phi_\epsilon}{\dee\bar z}(z-w)\ f(w)\ dA,$$
and the observation about the radially symmetric
function and our hypothesis about $f$ allows
us write the integral in polar coordinates
about $z$ to conclude that $f_\epsilon$ satisfies
the Cauchy-Riemann equations, and so is
analytic-diff on $D_{1-\epsilon}(0)$, and
consequently, is analytic-circ there, too. It
is easy to see that uniform limits of analytic-circ
functions are analytic-circ. We conclude that
$f$ is analytic-circ and so analytic in any
sense of the word.

\section{The Dirichlet problem in more general domains}
\label{sec7}
We solved the Dirichlet problem on the unit disc, given
polynomial boundary data, by explicitly extending
individual terms $z^n{\bar z}^{\,m}$ as harmonic polynomials.
Another way to approach this problem is via linear algebra.
Suppose a domain $\Om$ is described via a real polynomial
defining function $r(x,y)$ (meaning that
$\Om= \{x+iy\,:\, r(x,y)<0\}$), where $r(x,y)$ is 
of degree two. For $\Om$ to be a bounded domain,
it is clear that the boundary of $\Om$, which is the
zero set of $r$, must be a circle or an ellipse. Let
$\Delta$ denote the Laplace operator. Now,
the map $\mathcal F$ that takes a polynomial
$p(x,y)$ to the polynomial $\Delta(r p)$ maps the
finite dimensional vector space ${\mathcal P}_N$ of
polynomials of degree $N$ or less to itself. (Multiplying
by $r$ increases the degree by two, and applying the
second order operator $\Delta$ brings it back down by
two.) We claim that the map $\mathcal F$
is one-to-one on ${\mathcal P}_N$, and therefore
{\it onto}. Indeed, if $\Delta(r p)$
is the zero polynomial, then $r p$ is a harmonic
polynomial that vanishes on the boundary. The maximum
principle implies that it must be the zero polynomial.
Consequently, $p$ must be the zero polynomial, and
this proves that $\mathcal F$ is a one-to-one linear
mapping of a finite dimensional vector space into
itself, and so also onto. Now, to solve the Dirichlet
problem on $\Om$, given polynomial boundary data
$q(x,y)$, we know there is polynomial $p$ such that
$\Delta(r p)=\Delta q$. The polynomial $q-rp$ is
harmonic on $\Om$ and equal to $q$ on the boundary. It
solves the Dirichlet problem. We could have solved the
Dirichlet problem on the unit disc in the realm
of polynomials without ever writing a formula down!
Now we can solve the Dirichlet problem on an ellipse
using the same procedure that we did on the disc.
(However, the next obvious step, to try to write down
a Poisson integral formula on the ellipse, gets more
complicated because the monomials are not orthonormal
in the boundary inner product of the ellipse.)

The Khavinson-Shapiro conjecture states that discs and
ellipses are the only domains in the plane having
the property that solutions to the Dirichlet problem
with polynomial data must be polynomials. It is
tantalizing that it seems so much harder to settle
this question than the same problem with the word
``polynomial'' replaced by ``rational.'' Only discs
have the property that solutions to the Dirichlet
problem with rational boundary data must be rational
(see \cite{BEKS}).

For the remainder of this section, we will let
our Poisson and Dirichlet urges run rampant
and explain how the ideas in the previous sections
might be used to solve the Dirichlet problem on
more general domains. We will dispense with proofs
and follow a line of bold declarations. The
interested reader can find a more sober exposition
of some of these ideas in Chapters~22 and~34 of
\cite{Be} and in \cite{Be2}.

To solve the Dirichlet problem on a domain bounded
by a Jordan curve, one can use Carath\'eodory's
theorem (the theorem that states that the Riemann
map associated to such a domain extends continuously
to the boundary and maps the boundary one-to-one
onto the unit circle) to be able to pull back
solutions to the problem on the unit disc to the
domain. But we wonder if there might be a more
elementary way to do it.

Gustafsson's theorem \cite{G1} states that a bounded
finitely connected domain with $n$~continuous simple
closed boundary curves can be mapped to an $n$-connected
{\it quadrature domain\/} $\Om$ with smooth
real analytic boundary via a conformal mapping
that is continuous up to the boundary and as close to the
identity map in the uniform topology of the closure of
the domain as desired. Such a ``nearby'' quadrature domain
has the property that the average of an analytic function
over the domain with respect to area measure is a
finite linear combination of values of the function
and its derivatives at finitely many points in the
domain. The resulting ``quadrature identity'' is
the same for all analytic functions that are
square-integrable with respect to area measure on the
domain. Smooth real analytic curves have ``Schwarz
functions'' $S(z)$ that are analytic on a neighborhood
of the curve and satisfy $S(z)=\bar z$ on the curve.
The Schwarz functions associated to the boundary curves
of our quadrature domain $\Om$ have the following
stronger properties. There is a function $S(z)$ that
is meromorphic on an open set containing the closure of
$\Om$ that has no poles on the boundary of $\Om$
and satisfies the identity $\bar z=S(z)$ on the
boundary. Quadrature domains can be thought of as
a generalization of the unit disc (which is a
{\it one point\/} quadrature domain), and Gustafsson's
conformal mapping as a generalization of the Riemann
map in the $n$-connected setting.

Given a continuous function $\phi$ on the boundary
of our quadrature domain $\Om$, we can approximate
it by a rational function of $x$ and $y$ via the
Stone-Weierstrass theorem since the family of
such rational functions without singularities on
the boundary forms an algebra of continuous functions
that separates points. Writing such a rational
function as a rational function of $z$
and $\bar z$ and replacing $\bar z$ by $S(z)$
produces a meromorphic function on a neighborhood
of the closure of $\Om$ that has no poles on the
boundary and that agrees with the given rational
function on the boundary. If we can solve the
Dirichlet problem on $\Om$ with boundary data
$$\frac{1}{(z-a)^n}$$
for fixed $a$ in $\Om$ and positive integers $n$,
then, by subtracting such solutions from the
data, we would have harmonic functions that vanish
on the boundary and have general pole behavior at
$z=a$. We could then use these functions to subtract
off the poles of our meromorphic function and obtain
a solution to the Dirichlet problem with the given
rational boundary data. Then we could take uniform limits
and solve the Dirichlet problem for continuous
boundary data $\phi$ just like we did in the unit
disc.

In case $\Om$ is simply connected, it is possible
to solve the Dirichlet problem with boundary
data $(z-a)^{-n}$ using a Riemann mapping function.
Let $f:\Om\to D_1(0)$ be a Riemann map. The Green's
function $G(z,w)$ for $\Om$ is a constant times
$$\ln\left|\frac{f(z)-f(w)}
{1-\overline{f(w)}\,f(z)}\right|,$$
and derivatives
$$\frac{\dee^m}{\dee w^m}G(z,w)$$
are harmonic on $\Om-\{w\}$, are continuous
up to the boundary in $z$ and vanish on the
boundary in $z$, and the singularity at $w$
is precisely of the form a constant times
$(z-w)^{-m}$. One does not need to know that
the Riemann map is continuous up to the boundary
to see that these functions extend continuously
up to the boundary in $z$ and vanish there.
This follows from the fact that conformal
mappings are {\it proper\/} mappings: the
inverse image of a compact subset of the
unit disc is a compact subset of $\Om$.

Hence, in the simply connected case, we have
a method to solve the Dirichlet problem rather
analogous to the method we used in the case of
the unit disc. There is something appealing
about taking a close approximation to our
original domain, followed by a close
approximation to the boundary data, to be able
to find an elementary formula for the
solution to the Dirichlet problem.

Riemann maps associated to simply
connected quadrature domains can be expressed
as rational combinations of $z$ and the Schwarz
function, so solutions to the Dirichlet
problem with rational boundary data can also be
expressed as rational combinations of $z$ and
the Schwarz function!

Another way to construct the Poisson kernel is
to express it in terms of a normal derivative of
the Green's function, which, on a simply connected
quadrature domain, is also expressible in terms of
a Riemann map, and hence, also expressible in terms
of $z$ and the Schwarz function. It follows that the
Poisson kernel of a simply connected quadrature
domain is expressible in terms of $z$ and the
Scwharz function. Could we do similar things in
the multiply connected setting? Could we use
Ahlfors maps in place of a Riemann map? Might
the Poisson kernel there be expressible in terms
of $z$ and a Schwarz function and the harmonic
measure functions? We wonder.


\begin{thebibliography}{000}

\bibitem{Ah}
Ahlfors, L., {\em Complex analysis}, 3rd Edition, McGraw Hill,
New York, 1979.

\bibitem{Be}
Bell, S.,
{\em The Cauchy transform, potential theory, and
conformal mapping}, 2nd Edition, CRC Press, Boca Raton, 2016.

\bibitem{Be2}
Bell, S., {\em The Dirichlet and Neumann and Dirichlet-to-Neumann
problems in quadrature, double quadrature, and non-quadrature domains},
Analysis and Mathematical Physics {\bf 5} (2015), 113--135.

\bibitem{BEKS}
Bell, S., P.~Ebenfelt, D.~Khavinson, H.~Shapiro,
{\em On the classical Dirichlet problem in the plane with
rational data}, J. Anal. Math. {\bf 100} (2006), 157--190.

\bibitem{Gl}
Globevnik, J., {\em Zero integrals on circles and
characterizations of harmonic and analytic functions},
Trans. Amer. Math. Soc. {\bf 317} (1990), 313--330.

\bibitem{G1}
Gustafsson, B.,
{\em Quadrature domains and the Schottky double},
Acta Applicandae Math. {\bf 1} (1983), 209--240.

\bibitem{Ru}
Rudin, W., {\em Real and Complex Analysis}, McGraw Hill, New York, 1987.

\bibitem{Sp}
Springer, G., {\em On Morera's theorem}, Amer. Math. Monthly {\bf 64}
(1957), 323--331.

\bibitem{St}
Stein, E.~M. and R.~Shakarchi, {\em Complex analysis}, Princeton
Lectures in Analysis, Princeton Press, Princeton, 2003.

\end{thebibliography}
\end{document}